\author{\Large{Damanvir Singh Binner 
}}
\documentclass[12pt]{article}
\usepackage[margin=1in]{geometry}
\usepackage[usenames]{xcolor}
\usepackage{amssymb}
\usepackage{relsize}
\usepackage{amsmath}
\usepackage{amsthm}
\usepackage{amsfonts}
\usepackage{amscd}
\usepackage{graphicx}
\usepackage{hyperref}
\usepackage{tikz-cd}
\usepackage{tikz}

\begin{document}

\newcommand{\D}[1]{{\bf \color{red} #1}}
\newcommand{\M}[1]{{\bf \color{magenta} #1}}
\newcommand{\K}[1]{{\bf \color{violet} #1}}

\theoremstyle{plain}
\newtheorem{theorem}{Theorem}
\newtheorem{corollary}[theorem]{Corollary}
\newtheorem{lemma}[theorem]{Lemma}
\newtheorem{proposition}[theorem]{Proposition}
\newtheorem{question}[theorem]{Question}

\theoremstyle{definition}
\newtheorem{definition}[theorem]{Definition}
\newtheorem{example}[theorem]{Example}
\newtheorem{conjecture}[theorem]{Conjecture}

\theoremstyle{remark}
\newtheorem{remark}[theorem]{Remark}

\title{\Large{Number of partitions of $n$ into parts not divisible by $m$}}
\date{}
\maketitle
\begin{center}
\vspace*{-8mm}
\large{Department of Mathematics \\
Indian Institute of Science Education and Research (IISER) \\
Mohali, Punjab, India \\
damanvirbinnar@iisermohali.ac.in}
\end{center}

\begin{abstract}
In this note, we obtain a formula which leads to a practical and efficient method to calculate the number of partitions of $n$ into parts not divisible by $m$ for given natural numbers $n$ and $m$. Our formula is a generalization of Euler's recurrence for integer partitions which can be viewed as the case $m=1$ of our formula. Our approach primarily involves the principle of inclusion and exclusion. We also use our approach to obtain a natural combinatorial proof of a identity of Glaisher which generalizes a classical theorem of Euler.
\end{abstract}

\section{Introduction}
\label{S1}

Euler \cite[Corollary 1.2]{Andrews} proved that the number of partitions of $n$ into distinct parts is equal to the number of partitions into distinct parts. Glaisher \cite{Glaisher, Konan} generalized Euler's result by proving that for any positive integers $m$ and $n$, the number of partitions of $n$ into parts not divisible by $m$ is equal to the number of partitions of $n$ with each part appearing less than $m$ times. In this note, we obtain a formula that helps us to quickly calculate these numbers. Our formula generalizes Euler's recurrence relation for integer partitions \cite[Corollary 1.8]{Andrews}. We use the principle of inclusion and exclusion (PIE) to study these objects, an approach used by the present author and Rattan in the author's PhD Thesis \cite[Section 5.1]{Thesis} to obtain a natural combinatorial proof of Euler's recurrence. This approach was also used by the present author \cite{Parity} to obtain a formula for the number of partitions of $n$ with a given parity of the smallest part.

Note that the calculation of $p(n)$ is very easy compared to finding all the partitions of $n$ because of the availability of formulae such as Hardy-Ramanujan-Rademacher formula and recurrences such as Euler's recurrence. Therefore, Theorem \ref{Gen} below describes an efficient method to calculate the number of partitions of $n$ into parts not divisible by $m$, as demonstrated by examples in Section \ref{eg}.

\section{Main Theorem}
\label{main}

\begin{theorem} 
\label{Gen}
Let $P_m(n)$ denote the number of partitions of $n$ into parts not divisible by $m$. Then, $P_m(n)$ is given by the following formula.
$$ P_m(n) = p(n) + \sum_{k \geq 1} (-1)^k \left(p\left( n-\frac{mk(3k-1)}{2} \right) + p\left( n-\frac{mk(3k+1)}{2} \right) \right). $$
\end{theorem}

\begin{remark}
For $m=1$, $P_m(n)=0$, and Theorem \ref{Gen} immediately yields Euler's recurrence. 
\end{remark}

\begin{proof}[Proof of Theorem \ref{Gen}]

We recall some notation defined in \cite[Section 5.1]{Thesis} and define some new notation. 

\begin{itemize}
\item $A_{j,k}(n)$ is the set of partitions of $n$ having exactly $k$ parts of size $j$;
\item $B_{j,k}(n)$ is the set of partitions of $n$ having at least $k$ parts of size $j$. 
\item $C_{j,k}(n)$ is the set of partitions of $n$ having at most $k$ parts of size $j$;
\end{itemize}

The following properties of these sets are immediate. 

\begin{enumerate}
\item $|B_{j,k}(n)| = p(n-jk)$.
\item If $j \neq j'$, then $|B_{j,k}(n) \cap B_{j',k'}(n)| = p(n-jk-j'k')$.
\item $C_{j,k}(n) = B_{j,k+1}^{\complement}(n)$, where the complementation is with respect to the set $Par(n)$, consisting of all partitions of $n$. 
\item In particular, $A_{j,0}(n) = C_{j,0}(n) = B_{j,1}^{\complement}(n)$.
\end{enumerate}

We also need the following notation.

\begin{itemize}
\item D denotes the set of nonempty distinct partitions. 
\item For $s \in \mathbb{N}$, $T_s$ denotes the set of partitions into $s$ distinct parts.
\item For a partition $\pi$, $n(\pi)$ denotes the number of parts in $\pi$. 
\end{itemize}

Then using PIE, we have
\begin{align*}
 P_m(n) &= |A_{m,0}(n) \cap A_{2m,0}(n) \cap A_{3m,0}(n) \cap \cdots| \\
 &= |\cap_{i \geq 1} A_{im,0}(n)| \\
 &= |\cup_{i \geq 1} B_{im,1}^{\complement}(n)| \\
 &= \sum_{s \geq 0} (-1)^s \sum_{(i_1,i_2, \cdots, i_s) \in T_s} |B_{i_1m,1}(n) \cap B_{i_2m,1}(n) \cap \cdots \cap B_{i_sm,1}(n)| \\
 &= \sum_{s \geq 0} (-1)^s \sum_{(i_1,i_2, \cdots, i_s) \in T_s} p(n-i_1m-i_2m- \cdots i_sm) \\
 &= \sum_{s \geq 0} \sum_{\pi \in T_s} (-1)^s p(n-m|\pi|) \\
 &= p(n) + \sum_{s \geq 1} \sum_{\pi \in T_s} (-1)^s p(n-m|\pi|) \\
 &= p(n) + \sum_{\pi \in D} (-1)^{n(\pi)} p(n-m|\pi|),
 \end{align*}
 which completes the proof by Euler's pentagonal number theorem \cite[Theorem 1.6]{Andrews}. Note that here we are able to use PIE even though there are infinitely many sets because all except finitely many are empty. 
\end{proof}

\section{A proof of Glaisher's identity}
\label{Glaisher}

Let $Q_m(n)$ denote the number of partitions of $n$ with each part appearing less than $m$ times. Then using PIE, we have
\begin{align*}
 Q_m(n) &= |\cap_{i \geq 1} C_{i,m-1}(n)| \\
 &= |\cup_{i \geq 1} B_{i,m}^{\complement}(n)| \\
 &= \sum_{s \geq 0} (-1)^s \sum_{(i_1,i_2, \cdots, i_s) \in T_s} |B_{i_1,m}(n) \cap B_{i_2,m}(n) \cap \cdots \cap B_{i_s,m}(n)| \\
 &= \sum_{s \geq 0} (-1)^s \sum_{(i_1,i_2, \cdots, i_s) \in T_s} p(n-i_1m-i_2m- \cdots i_sm) \\
 &= \sum_{s \geq 0} \sum_{\pi \in T_s} (-1)^s p(n-m|\pi|) \\
 &= p(n) + \sum_{s \geq 1} \sum_{\pi \in T_s} (-1)^s p(n-m|\pi|) \\
 &= p(n) + \sum_{\pi \in D} (-1)^{n(\pi)} p(n-m|\pi|),
 \end{align*}
completing the proof of Glaisher's generalization of Euler's identity.

\section{Examples}
\label{eg}

For computational purposes, it is convenient to note a few terms and compare the pattern with the terms appearing in Euler's recurrence for partitions, which is given as
 \begin{multline*}
 p(n) = p(n-1)+p(n-2)-p(n-5)-p(n-7) \\
  + p(n-12) + p(n-15) - p(n-22) - p(n-26) + \cdots
  \end{multline*}
  
  Then, by Theorem \ref{Gen}, we have the following formula for $P_m(n)$.
 \begin{multline*}
 P_m(n) = p(n) - p(n-m) - p(n-2m) + p(n-5m) + p(n-7m) \\
  - p(n-12m) - p(n-15m) + p(n-22m) + p(n-26m) - \cdots
  \end{multline*}
  
  First suppose $n=17$ and $m=3$. Then, we have 
  \begin{align*}
   P_3(17) &= p(17) - p(14) - p(11) + p(2) \\
   &= 297 - 135 - 56 + 2 \\
   &= 108. 
   \end{align*}
   
   Thus, $108$ out of $297$ partitions of $17$ have parts not divisible by $3$ and the remaining $189$ partitions have at least one part that is divisible by $3$.
   
   Next, suppose $n=164$ and $m=7$. Then, we have 
  \begin{align*}
   P_7(164) &= p(164) - p(157) - p(150) + p(129)+p(115)-p(80)-p(59)+p(10) \\
   &= 156919475295 - 80630964769 - 40853235313 \\
   &+ 4835271870 + 1064144451 - 15796476 - 831820 + 42 \\
   &= 41318063280. 
   \end{align*}

Thus, $41318063280$ out of $156919475295$ partitions of $164$ have parts not divisible by $7$ and the remaining $115601412015$ partitions have at least one part that is divisible by $7$.

\end{document}